\input amstex
\documentstyle {amsppt}
\magnification=1200
\nologo
\vsize=8.9truein
\hsize=6.5truein
\nopagenumbers

\input xypic
\xyoption{all}
\xydashfont

\def\abs#1{\left\vert #1 \right\vert}

\topmatter

\title
The true prosoluble completion of a group
\\
Examples and open problems
\endtitle

\rightheadtext{The true prosoluble completion of a group}

\author
Goulnara Arzhantseva, Pierre de la Harpe, and Delaram Kahrobaei
\endauthor

\thanks
The authors acknowledge support from the
{\it Swiss National Science Foundation}.
\endthanks




\abstract
The {\it true prosoluble completion}
$P\Cal S (\Gamma)$ of a group $\Gamma$
is the inverse limit of the projective system of
soluble quotients of~$\Gamma$.
Our purpose is to describe examples
and to point out some natural open problems.
We discuss a question of Grothendieck
for profinite completions
and its analogue for true prosoluble
and true pronilpotent completions.

1. Introduction.

2. Completion with respect to a directed set of  normal subgroups.

3. Universal property.

4. Examples of directed sets of normal subgroups.

5. True prosoluble completions.

6. Examples.

7. On the true prosoluble and the true pronilpotent analogues
of Grothendieck's problem.


\endabstract

\address
Goulnara Arzhantseva and Pierre de la Harpe,
Section de Math\'ematiques, Universit\'e de Gen\`eve, C.P. 64,
CH--1211 Gen\`eve 4, Suisse.
\newline
E--mail: Goulnara.Arjantseva\@math.unige.ch,
Pierre.delaHarpe\@math.unige.ch
\endaddress

\address
Delaram Kahrobaei,
Mathematical Institute (Room 315), University of St Andrews,
North Haugh, St Andrews, Fife, KY16 9SS, Scotland.
\newline
E--mail: delaram.kahrobaei\@st-andrews.ac.uk
\endaddress

\endtopmatter

\document

\head{\bf
1.~Introduction
}\endhead

A group $\Gamma$ has a {\it profinite topology},
for which the set $\Cal F$ of normal subgroups of finite index
is a basis of neighbourhoods of the identity,
and the resulting {\it profinite completion},
hereafter denoted by $P \Cal F (\Gamma)$.
The canonical homomorphism
$$
\varphi_{\Cal F} : \Gamma \longrightarrow P \Cal F (\Gamma)
$$
is injective if the group $\Gamma$ is
{\it residually finite} (by definition).
The notion of profinite completion is relevant
in various domains (outside pure group theory),
including  Galois theory of infinite fields extensions
and fundamental groups in algebraic topo\-logy \cite{Groth--70}.
For the theory of profinite groups, there are
many papers and several books available
\cite{DiSMS--91}, \cite{RibZa--00}, \cite{Wilso--98};
see Section 1.1 in \cite{Se--73/94} for a quick introduction
and \cite{HallM--50} for an early paper.

Besides $\Cal F$,
there are other natural families of normal subgroups of $\Gamma$
which give rise to other \lq\lq procompletions\rq\rq .
The purpose of this report is to consider some variants,
with special emphasis on the {\it true prosoluble completion}
$P \Cal S (\Gamma)$ associated to the family $\Cal S$
of all normal subgroups of $\Gamma$ with soluble quotients.
The corresponding homomorphism
$$
\varphi_{\Cal S} : \Gamma \longrightarrow P \Cal S (\Gamma)
$$
is injective if the group $\Gamma$ is {\it residually soluble}
(by definition).

\medskip

On the one hand,
we discuss examples including free groups, free soluble groups, wreath
products, $SL_d(\bold Z)$ and its congruence subgroups, the  
Grigorchuk group,
and parafree groups.  
On the other hand,
we discuss some open problems,
of which we would like to point out from the start
the following  ones.

\medskip

(i) Let $\Gamma,\Delta$ be two residually finite groups and let
$\psi : \Gamma \longrightarrow \Delta$ be a homomorphism such that,
at the profinite level, the corresponding homomorphism
$P \Cal F (\psi) : P \Cal F (\Gamma) \longrightarrow P \Cal F (\Delta)$
is an isomorphism.
How far from an isomorphism can $\psi$ be?
The problem goes back to Grothendieck \cite{Groth--70},
and is  motivated by the need to compare
two notions of a fundamental group for algebraic varieties.
There are examples with $\psi$ not an isomorphism
and $P \Cal F (\psi)$ an isomorphism,
with $\Gamma,\Delta$ finitely generated~\cite{PlaTa--86}
and finitely presented \cite{BriGr--04}.
If $P \Cal F (\psi)$ is an isomorphism,
there are known sufficient conditions
\footnote{
Proposition 2 of \cite{PlaTa--90}:
if $\Gamma,\Delta$ are finitely generated residually finite groups,
if $\Delta$ is a soluble subgroup of $GL_n(\bold C)$ for some $n \ge 1$,
and if $P \Cal F (\psi)$ is an isomorphism,
then so is $\psi$.
}
on $\Gamma$ and $\Delta$
for $\psi$ to be an isomorphism.

Our main interest in this paper is  the
{\it prosoluble analogue of Grothendieck's problem.}
More precisely, let again $\Gamma,\Delta$ be two groups and
$\psi : \Gamma \longrightarrow \Delta$ a homomorphism,
but assume now that the groups are residually soluble.
Assume that, at the true prosoluble level,
the corresponding homomorphism
$P \Cal S (\psi) : P \Cal S (\Gamma) \longrightarrow P \Cal S (\Delta)$
is an isomorphism.
How far from an isomorphism can $\psi$ be?
Additional requirements can be  added
on $\Gamma$ and $\Delta$
(such as finite generation, finite presentation, ...).


\medskip

(ii) True prosoluble completions provide a natural setting to turn  
qualitative
statements  of the kind
\lq\lq some group $\Gamma$ is not residually soluble\rq\rq \
in more precise statements concerning the {\it true prosoluble kernel}
$\operatorname{Ker}(\varphi_{\Cal S} : \Gamma \longrightarrow P \Cal  
S (\Gamma))$.
One example is worked out in~(6.F).

\medskip

(iii) Can one find interesting characterizations of those groups
which are true prosoluble completions of residually soluble group?
More precisely, let $G$ be a complete Hausdorff topological group  
such that,
for any $g \in G$, $g \ne 1$,
there exists an open normal subgroup $N$ of $G$ not containing $g$
such that $G/N$ is soluble;
how can it be decided whether $G$ is isomorphic to $P \Cal S
(\Gamma)$ 
for some finitely generated group $\Gamma$?
for some finitely presented group $\Gamma$?
The corresponding questions for profinite groups are standard, and mostly open~\cite{KasNi--06}.

\medskip

Other open problems occur in (4.H), (6.F), (6.G), (7.B), and (7.C).

\medskip

We are grateful to Gilbert Baumslag, Slava Grigorchuk, Said Sidki,
and John Wilson  for valuable remarks. 
We also thank Dan Segal and an anonymous referee for pointing out
a mistake in a previous version of this paper.

\bigskip
\head{\bf
2.~Completion with respect to a directed set of  normal subgroups
}\endhead

In this section, we review some classical constructions and facts.
See in particular
\cite{Weil--40}, with \S~5 on projective limits,
\cite{Bourb--60}, with Chapter 3, \S~3, No.~4 on completions
and \S~7 on projective limits,
and \cite{Kelle--55}, with Problem~Q of Chapter~6 on completions.
Defining a topology on a group using a family of subgroups goes back  
at least to
Garrett Birkhoff \cite{Birkh--37, see Pages 52--54}
and Andr\'e Weil.

\medskip

Let $\Gamma$ be a group.
Let $\Cal N$ be a family of normal subgroups of $\Gamma$ which is  
directed,
namely which is such that the intersection of two groups in $\Cal N$
contains always a group in $\Cal N$.

\medskip

(2.A) Denote by $\Cal C \Cal N (\Gamma)$ the intersection of all  
elements in $\Cal N$
(the letter $\Cal C$ stands for \lq\lq core\rq\rq ),
by $\underline{\Gamma}$ the quotient  group
$\Gamma / \Cal C \Cal N (\Gamma)$,
and by $\underline{\Cal N}$ the family of  normal subgroups
of $\underline{\Gamma}$ which are images of groups in $\Cal N$.
Then $\underline{\Cal N}$ is a basis of neighbourhoods of the identity
for a Hausdorff topology on $\underline{\Gamma}$.
The corresponding left and right uniformities have the same Cauchy nets;
indeed, for $x,y \in \underline{\Gamma}$
and $\underline{N} \in \underline{\Cal N}$,
we have $x^{-1}y \in \underline{N}$ if and only if $xy^{-1} \in  
\underline{N}$.
It follows that $\underline{\Gamma}$ can be completed,
say with respect to the left uniformity,
to a Hausdorff complete
\footnote{
Recall that a topological group $G$ is {\it complete}
if both its left and right uniform structures
are complete uniform structures,
or equivalently if {\it one} of these structures
is a complete uniform structure;
see \cite{Bourb--60}, Chapter 3, \S~3.
Let $\Gamma$ be a topological group
and let $G$ denote its completion, as a topological space,
with respect to the left uniform structure;
a sufficient condition for $G$ to be
a completion of $\Gamma$ as a topological group
is that the left and right uniform structures on $\Gamma$
have the same Cauchy nets;
see Theorem 1 of No.~4 in the same book.
}
group which is called here the
{\it pro--$\Cal N$--completion of $\Gamma$}
and which is denoted by $P\Cal N (\Gamma)$.
The canonical homomorphism
$$
\varphi_{\Cal N} \, : \, \Gamma \longrightarrow P\Cal N (\Gamma)
$$
has kernel $\Cal C \Cal N (\Gamma)$ and image dense in $P\Cal N  
(\Gamma)$.
For $N \in \Cal N$, the projection
$\Gamma / \Cal C \Cal N (\Gamma) \longrightarrow \Gamma/N$
extends uniquely to a conti\-nuous homomorphism
$$
p_N \, : \, P\Cal N (\Gamma) \longrightarrow \Gamma/N
$$
which is onto.

\medskip

(2.B) Observe that the following properties are equivalent
\roster
\item"(i)" $\{1\} \in \Cal N$,
\item"(ii)" the topology defined by $\Cal N$  on $\Gamma$ is discrete,
\item"(iii)" $P \Cal N (\Gamma) = \Gamma$ and
$\varphi_{\Cal N} : \Gamma \longrightarrow P \Cal N (\Gamma)$
is the identity.
\endroster

More generally, $\Cal C \Cal N (\Gamma) \in \Cal N$ if and only if the topology
defined by $\underline{\Cal N}$ on $\Gamma / \Cal C \Cal N (\Gamma)$ is discrete, if and 
only if the natural homomorphism $\Gamma / \Cal C \Cal N (\Gamma)\longrightarrow P \Cal N (\Gamma)$ 
is the identity.

\medskip

(2.C)
Assume moreover that the family $\Cal N$ is countable.
We can assume without loss of generality that the elements of $\Cal N$
constitute a nested sequence
$N_1 \supset N_2 \supset \cdots$
(otherwise, if $\Cal N = \big\{ \tilde N_j \big\}_{j \ge 1}$,
choose inductively $N_{j+1}$ as a group in $\Cal N$
contained in $N_j \cap \tilde N_1 \cap \cdots \cap \tilde N_{j+1}$).
Assume also that $\cap_{i\ge1}N_i = \{1\}$.

The topology defined by $\Cal N$ on $\Gamma$ is metrisable.
Indeed, define first $v_i : \Gamma \longrightarrow \{0,1\}$
by $v_i(\gamma) = 0$ if $\gamma \in N_i$ and $v_i(\gamma) = 1$ if not.
Define next $w : \Gamma \longrightarrow [0,1]$
by $w(\gamma) = \sum_{i \ge 1}2^{-i}v_i(\gamma)$.
Then the mapping $d : \Gamma \times \Gamma \longrightarrow [0,1]$
defined by $d(\gamma_1,\gamma_2) = w(\gamma_1^{-1}\gamma_2)$
is a left--invariant ultrametric on $\Gamma$
which defines the same topology as $\Cal N$.

\medskip

(2.D) The quotients $$
\Gamma/N
\qquad \text{where $N \in \Cal N$,}
$$
and the canonical projections
$$
p_{M,N} \, : \, \Gamma/N \longrightarrow \Gamma/M
\qquad \text{where $M,N \in \Cal N$ are such that $N \subset M$,}
$$
constitute an inverse system of groups
of which the inverse limit (also called the projective limit)
$\varprojlim \Gamma/N$ can be identified with $P \Cal N (\Gamma)$.
The following properties are standard:
$P \Cal N (\Gamma)$ is totally disconnected and complete.

For a group $\Gamma$, denote by $Z(\Gamma)$ its center
and by $D(\Gamma)$ the subgroup generated by the commutators;
for a topological group $G$,  denote by  $\overline{D}(G)$ the  
closure of $D(G)$.
We have:
$$
Z(P \Cal N (\Gamma)) = \varprojlim Z(\Gamma/N)
\qquad \text{and} \qquad
\overline{D}(P \Cal N (\Gamma)) = \varprojlim D(\Gamma/N)
$$
(see \cite{Bourb--82}, Appendix I, No. 2).
See also (5.C) below.

\medskip

(2.E) Let $A$ be a partially ordered set which is directed.
Consider an inverse system consisting of groups $\Gamma_{\alpha}$,
with $\alpha \in A$, and homomorphisms
$p_{\alpha,\beta} : \Gamma_{\beta} \longrightarrow \Gamma_{\alpha}$,
with $\alpha,\beta \in A$ and $\alpha \le \beta$.
Let $\Gamma = \varprojlim \Gamma_{\alpha}$ denote the inverse limit.
Even when the groups $\Gamma_{\alpha}$ are not all reduced to $\{1\}$
and  the homomorphisms $p_{\alpha,\beta}$ are all onto,
the natural homomorphisms $p_{\alpha} : \Gamma \longrightarrow \Gamma_ 
{\alpha}$
need not be onto,
and indeed the limit $\Gamma$ can be reduced to one element;
see \cite{HigSt--54}.
(If $A$ is the set of natural integers, with the usual order,
it is straightforward to check that the $p_{\alpha}$ are onto.)

Because of this kind of phenomena, we have chosen here
to define pro--$\Cal N$--completions as appropriate topological
completions. However, it would be possible and equivalent to use inverse limits 
systematically.

\bigskip
\head{\bf
3.~Universal property
}\endhead

(3.A) Let $\Gamma$, $\Cal N$, and $P \Cal N (\Gamma)$ be as in the  
previous section.
Let $H$ be a topological group.
Assume that there is given a family
$\left( \psi_N : H \longrightarrow \Gamma/N \right)_{N \in \Cal N}$
of continuous homomorphisms such that,
for $M,N \in \Cal N$ with $N \subset M$,
the diagram
$$
\aligned
\xymatrix{ H \ar[rr]^- {\psi_N}  \ar[rd]^-{\psi_M} & & \Gamma /N
\ar[dl]_-{p_{M,N}} \\
& \Gamma /M }  
\endaligned
$$
commutes.
Then there exists a unique continuous homomorphism
$\psi : H \longrightarrow P \Cal N (\Gamma)$ such that the diagram
$$
\aligned
\xymatrix{ H \ar[rr]^- {\psi}  \ar[rd]^-{\psi_M} & & P\Cal N  (\Gamma )
\ar[dl]_-{p_{M}} \\
& \Gamma /M }
\endaligned
$$
commutes for all $M \in \Cal N$.

The pro--$\Cal N$--completion $P \Cal N (\Gamma)$ is characterized
up to unique continuous homomorphism
by this universal property.

\smallskip

{\it Caveat.} Even if $\psi_N$ is onto for each $N \in \Cal N$,
$\psi$ needs not be onto; see e.g. (6.A) below.
\medskip

 (3.B) Let $\Cal M$ be a subset of $\Cal N$
which is {\it cofinal}, namely such that any $N \in \Cal N$
contains some $M$ in $\Cal M$.
Then $\Cal C \Cal M (\Gamma) = \Cal C \Cal N (\Gamma)$,
the topology defined on $\Gamma / \Cal C \Cal M (\Gamma)$ by $\Cal M$
coincides with that defined by $\Cal N$,
and the continuous homomorphism
$P \Cal M (\Gamma) \longrightarrow P \Cal N (\Gamma)$
defined above is an isomorphism.
\medskip

    (3.C) Let $\Cal M$ [respectively $\Cal N$] be
a directed family of normal subgroups
of a  group $\Gamma$ [respectively $\Delta$]
and let $\psi : \Gamma \longrightarrow \Delta$
be a group homomorphism.
Assume that, for each $N \in \Cal N$,
the family 
$\Cal M_N=\left\{M\in \Cal M \mid \psi(M)\subset N\right\}$ is cofinal in $\Cal M$.

For a given $N \in \Cal N$,
the family of homomorphisms
$\Gamma/M \longrightarrow \Delta/N$ induced by $\psi$
(the family is indexed by $\Cal M_N$)
gives rise to a continuous homomorphism
$P \Cal M (\Gamma) \longrightarrow \Delta/N$.
In turn, these give rise to
a continuous homomorphism
$P \Cal M (\Gamma) \longrightarrow P \Cal N (\Delta)$.

This will occur several times in 
Section 4.

\medskip

\bigskip
\head{\bf
4.~Examples of directed sets of normal subgroups
}\endhead

(4.A) The set $\Cal F$ of normal subgroups of finite index
in a group $\Gamma$
gives rise to the {\it profinite completion}
$P \Cal F (\Gamma)$ of $\Gamma$.
There is a large literature on these completions,
alluded to in the introduction.

For a prime number $p$, the set $\Cal F_p$
of normal subgroups of index a power of $p$ in a group $\Gamma$
gives rise to the {\it pro--$p$--completion} $P_{\hat p}(\Gamma)$.
See \cite{DiSMS--91} and \cite{SaSeS--00}.
Since $\Cal F_p \subset \Cal F$, there is a canonical homomorphism
$$
P \Cal F (\Gamma) \, \longrightarrow \, P_{\hat p}(\Gamma)
$$
by (3.C). The resulting homomorphism
$P \Cal F (\Gamma) \longrightarrow \prod_p P_{\hat p}(\Gamma)$
is sometimes an isomorphism, as it is the case for $\bold Z$,
and sometimes not, as it is the case for
$\bigoplus_{n\ge 5}\operatorname{Alt}(n)$,
or for any non--trivial direct sum of non abelian finite simple groups.

\medskip

(4.B) The set $\Cal S$ of normal subgroups with soluble quotients
gives rise to the {\it true pro\-soluble completion} $P \Cal S  
(\Gamma)$ of $\Gamma$;
it is the main subject of the present note.

The {\it prosoluble completion} of the literature
refers usually to the family $\Cal F \Cal S$ of normal subgroups
with {\it finite} soluble quotients
(an exception is \cite{CocHa}, where our $P \Cal S (\Gamma)$
is called the \lq\lq pro--solvable completion\rq\rq \ of $\Gamma$,
and where other variants appear).
Since $\Cal F \Cal S \subset \Cal F$ and  $\Cal F \Cal S \subset \Cal  
S$,
there are canonical homomorphisms
$$
P \Cal F (\Gamma) \, \longrightarrow \, P \Cal F \Cal S (\Gamma)
\qquad \text{and} \qquad
P \Cal S (\Gamma) \, \longrightarrow \, P \Cal F \Cal S (\Gamma)
$$
by (3.C). See (6.A) and (6.B) for  examples.

The adjective \lq\lq true\rq\rq \ should not mislead the reader:
for some groups, for example for an infinite cyclic group,
the true prosoluble completion is \lq\lq much smaller\rq\rq \
than the prosoluble completion.

\medskip

(4.C) Similarly, we distinguish the
{\it true pronilpotent completion} $P \Cal Ni (\Gamma)$
from the pro\-nilpotent completion $P \Cal F \Cal Ni (\Gamma)$ of the  
literature.
There are again canonical homomorphisms
$$
P \Cal F (\Gamma)\, \longrightarrow\,
P \Cal F \Cal S (\Gamma) \,\longrightarrow\,
P \Cal F \Cal Ni (\Gamma)
\qquad \text{and} \qquad
P \Cal S (\Gamma) \, \longrightarrow \, P \Cal Ni (\Gamma)
\, \longrightarrow \, P \Cal F \Cal N i (\Gamma)
$$
by (3.C).

\medskip


%
%

\medskip

(4.D)
If $\Cal N$ is any of the classes appearing in
Examples (4.A) to (4.C) above (see also (4.J) below), $P \Cal N$ is a functor.
This means that any homomorphism $\psi : \Gamma \longrightarrow \Delta$
factors through
$\Gamma / \Cal C \Cal N (\Gamma) \longrightarrow \Delta / \Cal C \Cal  
N (\Delta)$
and then extends to a continuous homomorphism
$$
P \Cal N (\psi) \, : \, P \Cal N (\Gamma) \longrightarrow P \Cal N  
(\Delta) .
$$
As described in the introduction for $\Cal F$ and $\Cal S$, the
{\it pro--$\Cal N$--analogue of Grothendieck's problem}
is to find examples of
homomorphisms $\psi : \Gamma \longrightarrow \Delta$,
with $\Gamma,\Delta$ residually $\Cal N$
and possibly subjected to some extra conditions
(such as finite generation or finite presentation),
such that $\psi$ {\it is not} an isomorphism and
such that $P \Cal N (\psi)$ {\it is one}.

\medskip

(4.E)
In a locally compact group which is totally discontinuous,
any neighbourhood of the identity contains an open subgroup
(Corollary 1 in Chapter 3, \S~4, No. 6 of \cite{Bourb--60}).
It follows that a group which is profinite is residually finite.

Let $p$ be a prime number. Let $G$ be a pro--$p$--group,
namely a profinite group in which any open subgroup
is of index a power of $p$; assume
\footnote{
This hypothesis cannot be deleted.
Compare with \cite{Pet--73} or Example 4.2.13 in \cite{RibZa--00}, where an infinitely generated profinite 
group $G$ which is not isomorphic to $P \Cal F (G)$ is constructed.}
that $G$ is finitely generated
(we mean as a topological group,
i.e. there exists a finite subset of $G$
which generates a dense subgroup).
It is a theorem of Serre that any subgroup of finite index in $G$
is open (Theorem 4.2.5 in \cite{Wilso--98}).
It follows that, if $P_{\hat p}(G)$ denotes
the pro--$p$--completion of $G$ viewed as an abstract group,
the canonical homomorphism
$\varphi_{\hat p} : G \longrightarrow P_{\hat p}(G)$
is a continuous isomorphism
(see Proposition 1.1.2 in \cite{Wilso--98}).

Consider in particular a finitely generated group $\Gamma$
which is residually a finite $p$--group and the embedding
$\varphi_{\hat p} : \Gamma \longrightarrow P_{\hat p}(\Gamma)$
in its pro--$p$--completion. Then
$P_{\hat p}(\varphi_{\hat p}) : P_{\hat p}(\Gamma)
\longrightarrow P_{\hat p}(P_{\hat p}(\Gamma))$
is a continuous isomorphism, by which we will from no on
identify $P_{\hat p}(P_{\hat p}(\Gamma))$ with $P_{\hat p}(\Gamma)$.

\medskip

(4.F) By the recent solution, due to Nikolov and Segal,
of a conjecture of Serre, any subgroup of finite index
in a finitely generated profinite group is open \cite{NikSe--03, NikSe--06a, NikSe--06b}.
As a consequence, the conclusion of (4.E) carries over
from pro--$p$--completions to profinite completions;
more precisely:

Consider a finitely generated group $\Gamma$
which is residually finite and the embedding
$\varphi_{\Cal F} : \Gamma \longrightarrow P \Cal F (\Gamma)$
in its profinite completion. Then
$P \Cal F(\varphi_{\Cal F}) : P \Cal F (\Gamma)
\longrightarrow P \Cal F ( P \Cal F (\Gamma))$
is a continuous isomorphism, by which we will from now on
identify $P \Cal F (P \Cal F (\Gamma))$ with $P \Cal F (\Gamma)$.

\medskip

(4.G) Say that a subgroup $\Delta$ of a profinite group $G$
has the {\it congruence subgroup property}
if any normal subgroup $N$ of $\Delta$
is of the form $N = M\cap\Delta$
for some open normal subgroup $M$ of $G$.

Let $p$ be a prime number, $G$ a pro--$p$--group,
and $\Delta$ a dense subgroup of $G$.
Observe that $\Delta$ is residually a finite $p$--group
(since $G$ has this property).
Assume moreover that $G$ is finitely generated.
Let $\psi : \Delta \longrightarrow G$
denote the inclusion and let
$P_{\hat p}(\psi) : P_{\hat p}(\Delta) \longrightarrow G$
denote the homomorphism induced by $\psi$ on the pro--$p$--completions
(recall that  we identify $P_{\hat p}(G)$ with $G$).

On the one hand,
the homomorphism $P_{\hat p}(\psi)$ is onto.
Indeed, for any open normal subgroup $M$ of $G$,
the composition of $\psi$
with the projection $G \longrightarrow G/M$
is onto, by density of $\Delta$ in $G$.
On the other hand,
the following properties are equivalent:
\roster
\item"(i)"
$\Delta$ has the congruence subgroup property;
\item"(ii)"
the pro--$p$--topology and the topology induced by $G$ coincide on 
$\Delta$;
\item"(iii)"
the homomorphism
$P_{\hat p}(\psi) : P_{\hat p}(\Delta) \longrightarrow G$
is an isomorphism.
\endroster
Indeed, (ii) implies (iii) because,
for the topology induced by $G$, the completion of $\Delta$
coincides with its closure;
and (iii) obviously implies (ii).
As the topology of a topological group
is determined by the neighbourhoods of the identity,
(i) and (ii) are merely reformulations of each other.

{\it Examples.} Let $n\ge 2$ be an integer prime to $p$.
Recall that the multiplication by $n$ is invertible
in the group $\bold Z_p$ of $p$--adic integers.
Consider the situation where
$\Gamma = \bold Z$, $P_{\hat p}(\Gamma) = \bold Z_p$,
and $\Delta = \frac{1}{n}\bold Z$.
The inclusion $\psi : \bold Z \longrightarrow \frac{1}{n}\bold Z$,
which is not an isomorphism,
induces a continuous automorphism of $\bold Z_p$,
which is the multiplication by $\frac{1}{n}$.
{\it
This shows that the pro--$p$--analogue
of Grothendieck's problem has a straightforward solution.
}

Consider  an irrational $p$--adic integer $x$
and let now $\Delta$ be the subgroup of $\bold Z_p$
generated by $\bold Z$ and $x$.
Then $\Delta \approx \bold Z^2$,
so that $P_{\hat p}(\Delta)$
is isomorphic to $\bold Z_p \oplus \bold Z_p$
and the homomorphism $P_{\hat p}(\Delta) \longrightarrow \bold Z_p$
induced by the inclusion is not an isomorphism.

\medskip

(4.H) Let $G$ be a finitely generated profinite group,
and $\Delta$ a dense subgroup of $G$. 
Let $\psi : \Delta \longrightarrow G$
denote the inclusion. As in (4.G), the corresponding continuous homomorphism
$P\Cal F (\Delta) \longrightarrow G$
is always onto, and it is an isomorphism
if and only if $\Delta$ has the congruence subgroup property.

\medskip

{\it Open problem.} Let $\Gamma$ be the first Grigorchuk group;
it is an infinite group which is residually a finite $2$--group,
and all its proper quotients are finite $2$--groups.
This group was introduced in \cite{Grigo--80};
see also Chapter VIII in \cite{Harpe--00}.
It follows that $\Gamma$ embeds in its profinite completion,
which coincides with its pro--$2$--completion,
its true prosoluble completion, and its prosoluble completion.
[Though it is a digression from our main theme,
let us point out that $P_{\hat 2} (\Gamma)$
coincides moreover with  the closure of $\Gamma$
in the compact automorphism group
of the rooted dyadic tree on which $\Gamma$ acts in the usual way
(Theorem 9 in \cite{Grigo--00}).]


The problem is to find an element $g$ in the complement
of $\Gamma$ in $P \Cal F(\Gamma)$ such that the
subgroup $\Delta$ of $P \Cal F(\Gamma)$
generated by $\Gamma$ and $g$
has the congruence subgroup property.
This would provide one more example of a homomorphism
$\psi : \Gamma \longrightarrow \Delta$
which is not an isomorphism
and which would be such that the homomorphism
$P \Cal F (\psi)$ is an isomorphism,
namely one solution to the original Grothendieck problem
of a different nature than the existing ones
(from \cite{PlaTa--86}, \cite{PlaTa--89}, and \cite{BriGr--04}).

\medskip

(4.I) There are other cases than those appearing above
which are potentially interesting, of which we mention here two more.

\medskip

For a group $\Gamma$, the following two properties are equivalent:
\roster
\item"(i)" for any pair $N_1,N_2$ of normal subgroups not reduced to $ 
\{1\}$
    in $\Gamma$,
\item""    we have $N_1 \cap N_2 \ne \{1\}$;
\item"(ii)" there does {\it not} exist normal subgroups
    $N,N_1,N_2 \ne \{1\}$ of $\Gamma$
\item""    such that $N = N_1 \times N_2$;
\endroster
(we leave it as an exercise to the
reader to check this).
For $\Gamma$ with these properties, the family $\Cal N^{\ne 1}$
of all normal subgroups distinct from $\{1\}$ gives rise to the
{\it pronormal topology} on $\Gamma$; see \cite{GlSoS}.

\medskip

(4.J) The family $\Cal A$ of normal subgroups with amenable quotients
gives rise to the {\it proamenable completion} $P \Cal A (\Gamma)$ of  
$\Gamma$.
By (3.C), there are canonical homomorphisms from $P \Cal A (\Gamma)$
to $P \Cal F (\Gamma)$ and to $P \Cal S (\Gamma)$.
The related notion of residual amenability occurs for example in
\cite{Clair--99} and \cite{EleSz}.

\medskip

\bigskip
\head{\bf
5.~True prosoluble completions
}\endhead

(5.A) 
%
    An obvious obstruction to the residual solubility of a group $ 
\Gamma$
is the existence of a perfect subgroup not reduced to one element.
Any group $\Gamma$ contains a unique maximal perfect subgroup,
that we denote by $P_{\Gamma}$; this follows from the fact
that a subgroup generated by two perfect subgroups is itself perfect.
Observe that $P_{\Gamma}$ is contained
in the intersection $D^{\infty}(\Gamma)$
of all the groups in the derived series of $\Gamma$,
but the inclusion can be strict.
(In fact, $P_{\Gamma}$ is the intersection of all the groups
in the so--called transfinite derived series of $\Gamma$,
but this transfinite series can be as long, without repetition,
as the cardinality of $\Gamma$ allows;
see \cite{Malc'--49}.)

\medskip

(5.B) For a topological group $G$ and an integer $n \ge 0$,
we denote by $\overline{D^n}(G)$ the $n^{\text{th}}$ group of the
topological derived series, defined inductively by
$\overline{D^0}(G)  =  G$ and
$\overline{D^{n+1}}(G)  =  \overline{[\overline{D^n}(G),\overline{D^n} 
(G)]}$,
where $\overline{[A,B]}$ stands for the closure of the subgroup of $G$
generated by commutators $a^{-1}b^{-1}ab$ with $a \in A$ and $b \in B$.
We denote  by $G_{(n)}$ the quotient $G/\overline{D^n}(G)$
and by $\sigma_{(n)} : G \longrightarrow G_{(n)}$ the canonical  
projection.

Recall that a topological group $G$ is \lq\lq topologically soluble\rq 
\rq ,
namely such that $\overline{D^n}(G) = \{1\}$ for $n$ large enough,
if and only if it is soluble, namely such that the $n^{\text{th}}$ term
$D^n(G)$ of its ordinary derived series is reduced to $\{1\}$ for $n$  
large enough
(see for example Chapter~III, \S~9, No.~1 in \cite{Bourb--72}).

\medskip

(5.C) For a group $\Gamma$ and an integer $n \ge 0$,
we claim that $P \Cal S (\Gamma)_{(n)}$ is canonically isomorphic to $ 
\Gamma_{(n)}$.
This follows from contemplation of the commutative diagram

$$
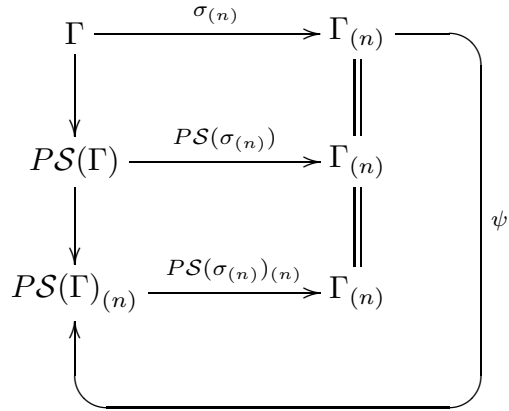
\aligned
\xymatrix{ \Gamma \ar[rr]^{\sigma_{(n)}}  \ar[d]   &  & \Gamma_{(n)} \ar@{=}[d] \ar `[r] `d[ddd]^{\psi} `l[lld] [ddll] &   \\
P\Cal S (\Gamma) \ar[rr]^-{P\Cal S(\sigma_{(n)})} \ar[d]      &  & \Gamma_{(n)} \ar@{=}[d]  & \\
{P\Cal S (\Gamma)}_{(n)}\ar[rr]^-{P\Cal S(\sigma_{(n)})_{(n)}} & & \Gamma_{(n)} & \\
& & & }
\endaligned
$$

\medskip

\noindent
where the existence and canonicity of $\psi$ follow from the  
universal property of
the projection $\sigma_{(n)}$.
In other words, as the morphism $\Gamma \longrightarrow P \Cal S  
(\Gamma)_{(n)}$
has a range which is soluble of degree $n$,
this morphism factors through $\Gamma_{(n)}$.

\medskip

(5.D) In a group $\Gamma$,
the family $\Cal S$ of all normal subgroups with soluble quotients
and the countable family $\left( D^n(\Gamma) \right)_{n \ge 0}$
define the same procompletion $P \Cal S (\Gamma)$, by Item (3.B).

    It follows from (2.C) that the topology on $P \Cal S (\Gamma)$
can be defined by a metric.

    Recall that the family $\Cal S$ can be uncountable.
P.~Hall has shown that this is the case for $\Gamma$  a free group of  
rank $2$;
see \cite{HallP--50}, and the exposition in \cite{BruBW--79}.

\medskip

(5.E) For a topological group $G$ and an integer $j \ge 0$,
we denote by $\overline{C^j}(G)$ the $j^{\operatorname{th}}$ group
of the topological lower central series,
defined inductively by  $\overline{C^1}(G)  =  G$ and
$\overline{C^{j+1}}(G)  =  \overline{[G,\overline{C^j}(G)]}$
for $j \ge 1$.
Then $G$ is nilpotent if and only if $\overline{C^j}(G) = \{1\}$
for $j$ large enough.
(Compare with (5.B).)

    The nilpotent quotients
$P \Cal N i (\Gamma)/\overline{C^j}(P \Cal N i (\Gamma))$
and $\Gamma / C^j(\Gamma)$ are isomorphic for all $j \ge 1$.
(Compare with (5.C).)

    The topology on $P \Cal N i (\Gamma)$ can be defined by a metric.
(Compare with  (5.D).)

\bigskip
\head{\bf
6.~Examples
}\endhead

(6.A) Let us show that the canonical homomorphism
$P \Cal S (\Gamma) \longrightarrow P \Cal F \Cal S (\Gamma)$
needs not be onto.

First, consider  an infinite cyclic group:  $\Gamma = \bold Z$.
Since $\bold Z$ is soluble,
$\varphi_{\Cal S} : \bold Z \longrightarrow P \Cal S (\bold Z)$
is an isomorphism onto.
As any finite quotient of $\bold Z$ is soluble (indeed abelian),
the prosoluble completion of $\bold Z$ coincides with its profinite  
completion.
Hence
$$
\bold Z \, = \, P \Cal S (\bold Z)
\, \not\approx \,
P \Cal F \Cal S  (\bold Z) \, \approx \,
P \Cal F  (\bold Z) \, \approx \,
\prod_p \bold Z_p.
$$

    Here is another family of examples.
For $k \ge 2$ and $d \ge 1$,
denote by $F(k,d)$ the quotient of a non--abelian free group
on $k$ generators by the $d$th term of its derived series
(the {\it free soluble group} of class $d$ with $k$ generators).
This group is obviously soluble and infinite,
and it is known to be residually finite (Theorem 6.3 in \cite 
{Gruen--57}).
Hence $F(k,d) =  P \Cal S (F(k,d))$ embeds properly in
its profinite completion $P \Cal F (F(k,d))$, and the latter  
coincides with
$P \Cal F \Cal S (F(k,d))$.

\medskip

(6.B) Similarly, the canonical homomorphisms
$P \Cal S (\Gamma) \longrightarrow P \Cal F \Cal S (\Gamma)$
needs not be injective.

Consider a wreath product $\Gamma = S \wr T$
where $S$ is soluble non--abelian and where $T$ is soluble infinite.
Theorem 3.1 in \cite{Gruen--57} establishes that $\Gamma$ is not  
residually finite,
so that in particular the morphism
$\varphi_{\Cal F \Cal S} : \Gamma \longrightarrow P \Cal F \Cal S  
(\Gamma)$
has kernel not reduced to $\{1\}$.
But $\Gamma$ is soluble, and therefore isomorphic to $P \Cal S  
(\Gamma)$.

For a finitely generated group $\Gamma$,
it is a theorem of P.~Hall  that $\Gamma / D^2(\Gamma)$
is residually finite \cite{HallP--59}, so that we have an embedding
$$
\Gamma / D^2(\Gamma) \, \approx \,
P \Cal S (\Gamma) / \overline{D^2}( P \Cal S (\Gamma) ) \,  
\longrightarrow \,
P \Cal F \left( \Gamma / D^2(\Gamma) \right) .
$$
In particular, a finitely generated group $\Gamma$ which is soluble  
of class $2$
always embeds in $P \Cal F (\Gamma)$.

However, there exists a finitely generated
soluble group of class $3$
which is non--Hopfian \cite{HallP--61}.
In particular, a finitely generated group $\Gamma$ which is soluble  
of class at least $3$
does not always embed in $P \Cal F (\Gamma) = P \Cal F \Cal S (\Gamma)$.
A finitely {\it presented} soluble group which is non--Hopfian is
described in \cite{Abels--79}.

\medskip

(6.C) For each integer $k \ge 2$,
let $F_k$ denote the nonabelian free group on $k$ generators.
This is a residually soluble group
(indeed it is residually a finite $p$--group for any prime $p$,
by a result of Iwasawa, see No.~6.1.9 in \cite{Robin--82}),
and therefore embeds in its true prosoluble completion $P \Cal S (F_k)$.

Since the abelianized groups $F_k/[F_k,F_k] \approx \bold Z^k$
are pairwise non--isomorphic, the true prosoluble completions
$P \Cal S (F_k)$ are pairwise non--isomorphic by (2.D) or (5.C).

\medskip

(6.D) For $d \ge 3$, the group $SL_d(\bold Z)$ is perfect
(this is an immediate consequence of the following fact:
any matrix in $SL_d(\bold Z)$ is a product of elementary matrices;
see for example Theorem 22.4 in \cite{MacDu--46}).
It follows that $P \Cal S (SL_d(\bold Z))$ is reduced to one element
and that $P \Cal S (GL_d(\bold Z))$ is the group of order two.

Consider however an integer $d \ge 2$,
a prime $p$, and the congruence subgroup
$$
\Gamma_d(p) \, = \, \operatorname{Ker}\left(
SL_d(\bold Z) \longrightarrow SL_d(\bold Z /p\bold Z)
\right),
$$
which is a subgroup of finite index in $SL_d(\bold Z)$.
Then $\Gamma_d(p)$ is residually  a finite $p$--group,
and therefore  embeds in $P \Cal S (\Gamma_d(p))$;
see for example Proposition 3.3.15 in \cite{RibZa--00},
since their proof written for $d=2$ carries over to any $d \ge 2$.

In particular, the property for a group $\Gamma$
to embed up to finite kernel in its prosoluble completion
is {\it not} stable by finite index.

\medskip



\medskip

(6.E)
Groups which are residually soluble and not soluble include
Baumslag--Solitar groups $\langle a,t \mid ta^pt^{-1} = a^q \rangle$
for $\abs{p},\abs{q} \ge 2$ \cite{RapVa--89},
positive one--relator groups \cite{Baums--71},
non--trivial free products
\footnote{
Recall of a particular case of a result of Gruenberg,
see Section~4 in \cite{Gruen--57}:
any free product of residually soluble groups is residually soluble.
By \lq\lq non--trivial\rq\rq \ free product, we mean that
the free product has at least two factors,
and that it is not the free product of two groups of order two.
}
of soluble groups,
and various just non-soluble groups \cite{BruSV--99}.
Also, some free products of soluble groups amalgamated over a cyclic  
group
are residually soluble and not soluble \cite{Kahro1}.

These are potential examples for further investigation
of the properties of true prosoluble completions.

\medskip

(6.F) For groups which are not residually soluble,
it is a natural problem to find out
properties of the true prosoluble kernel
$$
\operatorname{Ker}(\Gamma \longrightarrow P \Cal S (\Gamma))
\, = \,
\bigcap_{n=0}^{\infty} D^n(\Gamma) .
$$

For example, consider  the one--relator group
\footnote{
Recall that $[x,y] = x^{-1}y^{-1}xy$ and $x^y = y^{-1}xy$.
}
$$
\Gamma \, = \, \langle a, b ; a = [a,a^b] \rangle
$$
which appears in \cite{Baums--69} and \cite{Baums--71}.
Baumslag has shown that the profinite kernel of $\Gamma$
coincides with the derived group $D(\Gamma)$,
which is also the smallest normal subgroup of $\Gamma$ containing $a$,
and which is not reduced to one element by the Magnus Freiheitssatz.
As $D(\Gamma)$ is perfect, the true prosoluble kernel of $\Gamma$
coincides also with $D(\Gamma)$.
Since $\Gamma / D(\Gamma) \approx \bold Z$,
we have
$\bold Z \approx P \Cal S (\Gamma) \not\approx P \Cal F (\Gamma)
\approx \prod_p \bold Z_p$;
compare with (6.A).

\medskip

Examples to investigate include:
\roster
\item"$\circ$"
other one--relator groups which are not residually soluble;
\item"$\circ$"
wreath products as in (6.B);
\item"$\circ$"
free products of nilpotent groups with amalgamation
(see for example Proposition~7 in \cite{Kahro1});
\item"$\circ$"
parafree groups (see (7.B));
\item"$\circ$"
meta--residually soluble groups,
namely groups having a residually soluble normal subgroup
with residually soluble quotient
(see Theorem 2 in \cite{Kahro2});
particular cases: free--by--free groups which are not residually  
soluble.
\endroster

\medskip

(6.G) Does there exist a discrete group $\Gamma$
which on the one hand is resi\-dually finite and residually soluble,
namely for which both homomorphisms
$\Gamma \longrightarrow P \Cal F (\Gamma)$,
$\Gamma \longrightarrow P \Cal S (\Gamma)$
are injective,
and which on the other hand is such that the homomorphism
$\Gamma \longrightarrow P \Cal F \Cal S (\Gamma)$
is NOT injective?

\bigskip
\head{\bf
7.~On the true prosoluble and the true pronilpotent analogues
of Grothendieck's problem
}\endhead

(7.A)
Consider as in (4.H) the Grigorchuk group $\Gamma$
and  its pro--$2$--completion $P_{\hat 2}(\Gamma)$.

\proclaim{Hypothesis}  We assume from now on
      that there exists a finitely generated subgroup $\Delta$
      of $P_{\hat 2}(\Gamma)$ which contains $\Gamma$ properly
      and which has the congruence subgroup property.
\endproclaim

Observe that any finite quotient of $\Delta$ is a $2$--group,
as a consequence of the congruence subgroup property;
thus, the canonical homomorphism from $P \Cal F (\Delta)$
to $P_{\hat 2}(\Delta)$ is an isomorphism.
Let $\psi : \Gamma \longrightarrow \Delta$ denote the inclusion.
We know from (4.G) that
$P_{\hat 2}(\psi) : P_{\hat 2}(\Delta) \longrightarrow P_{\hat 
2}(\Gamma)$
is an isomorphism.

\proclaim{Lemma} With the notation above,
a quotient group of $\Delta$ is finite if and only if it is soluble.
\endproclaim

\demo{Proof}
Let $N$ be a normal subgroup of $\Delta$ such that $\Delta/N$ is finite.
We have already obseerved that
the quotient $\Delta/N$ is  a finite $2$--group,
in particular a nilpotent group, a fortiori a soluble group.

For the converse, we proceed by contradiction and assume that
there exists a normal subgroup $N$ in $\Delta$ such that the quotient
$\Delta /N$ is soluble and infinite.
Since $\Delta$ is finitely generated,
Zorn's lemma implies that there exists
a normal subgroup $M$ in $\Delta$ containing $N$
such that $\Delta/M$ is just infinite.

    Let $\Delta/M = D^0 \supset D^1 \supset D^2 \supset \cdots$
denote the derived series of the just infinite soluble group $\Delta/M$;
denote by $k$ the smallest integer
such that $D^{k+1}$ is of infinite index in  $\Delta/M$.
Then $D^{k+1} = \{1\}$ because $\Delta/M$ is just infinite.
The group $D^k$ is abelian and of finite index in $\Delta/M$;
it is therefore finitely generated.
The torsion subgroup of $D^k$ is normal in $\Delta/M$;
it is of infinite index in $D^k$, and thus also in $\Delta/M$;
hence $D^k$ has no torsion.
Let $\Delta_1$ denote the inverse image of $D^k$ in $\Delta$;
it is a normal subgroup of finite index containing~$M$.
The quotient $\Delta_1/M$
is finitely generated and free abelian,
say  $\Delta_1/M \approx \bold Z^d$ for some $d \ge 1$.

    Let $Q_2$ be the subgroup of $\Delta_1/M$ generated by the cubes,
and denote by $\Delta_2$ its inverse image in $\Delta_1$;
observe that $\Delta_2$ is of finite index in $\Delta_1$,
and therefore also in $\Delta$,
and that $\Delta_1 / \Delta_2 \approx (\bold Z /3\bold Z)^d$.
Let $\Delta_3$ be the intersection of all the conjugates of $\Delta_2$
in $\Delta$;
observe that $\Delta_3$ is of finite index in $\Delta_2$,
and therefore also in $\Delta$, and that it is also normal in $\Delta$.
Then $\Delta/\Delta_3$ is a finite quotient of $\Delta$
of order
$$
[\Delta : \Delta_1] \times [\Delta_1 : \Delta_2] \times [\Delta_2 : 
\Delta_3]
\, = \,
[\Delta : \Delta_1] \times 3^d \times [\Delta_2 : \Delta_3] .
$$
In particular, $3$ divides the order of $\Delta/\Delta_3$.

We have obtained a contradiction because,
since $P_{\hat 2}(\Delta) \approx P \Cal F (\Delta)$,
any finite quotient of $\Delta$ is a $2$--group.
This ends the proof.
\hfill $\square$
\enddemo

\proclaim{Consequence} Let $\Gamma$ be the Grigorchuk group.
If there would exist a dense subgroup
$\psi : \Delta \longrightarrow P_{\hat 2}(\Gamma)$
containing properly $\Gamma$ and
with the congruence subgroup property,
then the morphism
$$
P \Cal S (\psi) \, : \,
P \Cal S (\Gamma)  \, \longrightarrow \, P \Cal S (\Delta)
$$
induced by $\psi$ on the true prosoluble completions
would be an isomorphism.
\endproclaim

\medskip


\medskip

(7.B) Let us describe how Problem (i) from the introduction
has been studied for true pronilpotent  completions.

Recall that a group $\Gamma$ is {\it parafree}
if it is residually nilpotent and if there exists a free group $F$
such that   $F/C^j(F)$ and $\Gamma/C^j(\Gamma)$
are isomorphic for all $j \ge 1$.
Nonfree parafree groups have been discovered by G.~Baumslag
\cite{Baums--67};
later papers include \cite{Baums--68} and \cite{Baums--05}.

Let $\Gamma$ be a parafree group with finitely generated abelianization.
Let $F$ be as above; observe that $F$ is finitely generated.
Choose a subset $T$ of $\Gamma$ of which the canonical image
freely generates the free abelian group $\Gamma / C^2(\Gamma)$;
observe that $T$ is finite.
(Be careful: $T$ need not generate $\Gamma$.)
Let $S$ be a free set of generators of $F$
such that the canonical image of $S$ and $T$
are in bijection with each other
through the given isomorphism $F/C^2(F) \approx \Gamma/C^2(\Gamma)$,
and let $\varphi : S \longmapsto T$ be a compatible bijection.
Then $\varphi$ extends to a homomorphism, again denoted by $\varphi$,
from $F$ to $\Gamma$, and this $\varphi$ induces the given isomorphism
from $F / C^2(F)$ onto $\Gamma / C^2(\Gamma)$.
A group homomorphism with range a nilpotent group $A$ is onto
if and only if its composition with the abelianization
$A \longrightarrow A/C^2(A)$ is onto
(see e.g. \cite{Bourb--70}, Corollary 4, Page A I.70);
it follows that the homomorphism $\varphi_{(j)}$
from $F / C^j(F)$ to $\Gamma / C^j(\Gamma)$ induced by $\varphi$
is onto for all $j \ge 1$.
Since the group $F/C^j(F) \approx \Gamma/C^j(\Gamma)$
is Hopfian for all $j \ge 1$
(any finitely generated residually finite group is Hopfian,
by Mal'cev theorem \cite{Mal'c--40}),
it follows that $\varphi_{(j)}$ is an isomorphism for all $j \ge 1$.
We have shown:

\proclaim{Observation}
A  residually nilpotent group $\Gamma$ with finitely generated  
abelianization
is para\-free if and only if
there exist a free group $F$ of finite rank and a homomorphism
$\varphi : F \longrightarrow \Gamma$
which induces an isomorphism
$P \Cal N i (\varphi) : P \Cal N i (F) \longrightarrow P \Cal N i 
(\Gamma)$
on the true pronilpotent completions.
\endproclaim

Note that, in case the set $T$ is moreover generating for $\Gamma$,
the group $\Gamma$ itself is free;
see Problem 2 on Pages 346--347 of \cite{MaKaS--66}.
However, as G.~Baumslag discovered,
there are pairs $(\Gamma,F)$ as in the observation
with $\Gamma$ not free and generated by $k+1$ elements,
and with $F$ free of rank $k$,
for each $k \ge 2$.

    There is a related example on Page 173 of \cite{Stall--65}.
Let $F_2$ be the free group on two generators $x$ and $y$.
Set $y' = yxyx^{-1}y^{-1}$.
Let $\Gamma$ be the subgroup of $F_2$ generated by $x$ and $y'$.
Then $\Gamma$ is a proper subgroup of $F_2$, because $y \notin \Gamma$
(though $\Gamma$ is isomorphic to $F_2$).
The inclusion of $\Gamma$ in $F_2$ provides an isomorphism
$\Gamma/C^j(\Gamma) \longrightarrow F_2/C^j(F_2)$ for all $j \ge 1$,
and therefore an isomorphism from $P \Cal N i (\Gamma)$ onto $P \Cal  
N i (F_2)$.

\medskip

(7.C) In \cite{Baums--68}, there are examples of pairs $(F,\Gamma)$ of
groups with the following properties: $F$ is finitely generated and  
free,
both $F$ and $\Gamma$ are residually soluble,
the quotients $F/D^j(F)$ and $\Gamma/D^j(\Gamma)$ are isomorphic
for all $j \ge 0$, nevertheless $F$ and $\Gamma$ are {\it not}  
isomorphic
(indeed $\Gamma$ is not finitely generated).

However, nothing like the observation of Item (7.C) holds for
the quotients by the groups of the derived series.
We do not know if there exist
a finitely generated free group $F$,
a group $\Gamma$ which is residually soluble and not free,
and a homomorphism $\psi : F \longrightarrow \Gamma$
such that $\psi$ induces an isomorphism 
$F/D^j(F)   \longrightarrow   \Gamma/D^j(\Gamma)$
for all $j \ge 1$,
or equivalently such that
$P \Cal S (\psi) : P \Cal S (F) \longrightarrow P \Cal S (\Gamma)$
is an isomorphism.


\enddocument


\head{\bf
Groupes resolubles libres
}\endhead

D'apr\`es Geba, $F(k,d)$, $d \ge 3$,
n'est pas de pr\'esentation finie.
Voir A.L. Shmel'kin, {\it Wreath products and varieties of groups,}
Izv. Adak. Nauk. SSSR Ser. Mat. {\bf 29} (1965) 149--170.

\head{\bf
Extension du groupe de Grigorchuk
}\endhead

A propos de $\Gamma \ltimes \bold Z /2\bold Z$, voir l'article de
Said et Slava in IJAC. Ces groupes ne sont pas isomorhes, par exemple  
parce
que $(\Gamma \ltimes \bold Z /2\bold Z)_{\operatorname{ab}} \ne
(\Gamma)_{\operatorname{ab}}$ (un facteur $\bold Z /2\bold Z$ de plus).

\head{\bf
Anciennement partie de (7.B)
}\endhead

Let $\delta$ be an element of infinite order in $P_{\hat 2} (\Gamma)$.
Such elements exist, from direct inspection,
or as a consequence of Theorem 10 in \cite{Grigo--00}.
Let $\Delta$ be the subgroup of $P_{\hat 2}(\Gamma)$
generated by $\Gamma$ and $\delta$.
Observe that $\Delta$ is residually a finite $2$--group,
since it is a subgroup of $P_{\hat 2}(\Gamma)$.
It follows as above (see (4.F)) that the inclusion
of $\Delta$ in $P_{\hat 2}(\Delta)$ induces the identity automorphism
of $P_{\hat 2}(\Delta)$.

On the one hand,  the inclusion of $\Gamma$ in $\Delta$
is not an isomorphism since $\Gamma$ is a torsion group.
On the other hand, the sequence of inclusions
$$
\Gamma \, \subset \, \Delta  \, \subset \,
P_{\hat 2}(\Gamma)  \, \subset \, P_{\hat 2}(\Delta)
$$
provides for pro--$2$--completions a sequence of morphisms
$$
P_{\hat 2}(\Gamma)  \, \longrightarrow \, P_{\hat 2}(\Delta) \,  
\longrightarrow \,
P_{\hat 2}(\Gamma)  \, \longrightarrow \, P_{\hat 2}(\Delta)
$$
where the composition
$P_{\hat 2}(\Gamma)  \, \longrightarrow \, P_{\hat 2}(\Gamma)$
is the identity, and similarly for the composition
$P_{\hat 2}(\Delta)  \, \longrightarrow \, P_{\hat 2}(\Delta)$.
Thus the homomorphism
$P_{\hat 2}(\Gamma)  \, \longrightarrow \, P_{\hat 2}(\Delta)$
is also an isomorphism,
and this is precisely $P_{\hat 2}(\psi)$
if $\psi$ stands for the inclusion of $\Gamma$ in $\Delta$.

\bigskip
\head{\bf
Bribes non concluantes pour un (7.D) -- September 21st, 2005
}\endhead

(7.D) We use the same notation as in (7.B).
there exist an example of a nonfree parafree group $\Gamma$,
a free group $F$,
and a homomorphism $\varphi : F \longrightarrow \Gamma$
which is not an isomorphism and which induces an isomorphism
$P \Cal S (\varphi) : P \Cal S (F) \longrightarrow P \Cal S (\Gamma)$
on the true prosoluble completions.

The example of \cite{Baums--68} is not finitely generated;
we do not know any finitely generated example.

To end our Note, let us describe a very particular case of this example,
without reference to the vocabulary of varieties of groups used in  
\cite{Baums--68}.

Let $F$ be a free group on $k \ge 2$ generators $s_1,\hdots,s_k$.
For each integer $j \ge 1$, choose an element $\gamma_j \in D^j(F)$.
Denote by $\alpha_j$ the homomorphism from a copy $F_{j-1}$ to a copy  
$F_{j}$ of $F$
which maps $s_1$ to $s_1\gamma_j$ and $s_l$ to $s_l$ for $l \in \{2, 
\hdots,k\}$.
Observe that the image of $\alpha_j$ is a subgroup of $F_{j}$ which  
is free of rank $k$;
indeed, it is of rank at most $k$, because it is generated by $k$  
elements,
and of rank at least $k$, because its image in the abelian group $F_j/ 
D(F_j)$ is
of rank $k$.
Let $\Gamma$ be the direct limit of the inductive system defined by  
the $\alpha_j$.
Observe that $\Gamma$ is not finitely generated
(since any finitely generated subgroup of $\Gamma$ is contained in
$F_N$ for $N$ large enough, and this is a proper subgroup),
that it is locally free,
and that it is not free
(since $\Gamma/C^2(\Gamma) \approx F/C^2(F)$ is finitely generated).
We leave it to the reader to check that $\Gamma$ is a parafree group
and that $\Gamma / D^j(\Gamma) \approx F / D^j(F)$ for all $j \ge 1$.

\noindent
??????????????
Pas encore clair pour moi \`a l'heure o\`u j'\'ecris~:

s'il existe un homomorphisme $\beta$ de $F$ dans $\Gamma$ qui induise  
des isomorphismes
$ F / D^j(F) \approx \Gamma / D^j(\Gamma)$ pour tout $j \ge 1$.

%

\bigskip
\head{\bf
Notes -- September 7th, 2005
}\endhead

\subhead
(2.E)
\endsubhead

A propos de \cite{HigSt--54}, voir aussi Bourbaki, Th\'eorie des  
ensembles, Hermann 1970,
Page EIII.58.

Topologie sur un groupe d\'efinie par une famille de sous-groupes~:  
voir Bourbaki,
Topologie g\'en\'erale, Hermann, 1960, chapitre 3, \S~1, no 2, exemple.

Compl\'etude des limites projectives~: voir Bourbaki, Topologie g\'en 
\'erale, Hermann,
1965, chapitre 2, \S~3, no 5.

Compl\'etion de $\Gamma$ pour la $\Cal N$-topologie isomorphe \`a la  
limite projective~:
voir Bourbaki, Topologie g\'en\'erale, chapitre 3, \S~7, no 3,  
proposition 2.
(Et corollaire 2...)

\subhead
Anciennement en (5.A)
\endsubhead

    Digression
\footnote{
This is motivated by an example on Page 564 of
I.~Kapovich and D.~Wise,
{\it The equivalence of some residual properties of word-hyperbolic  
groups,}
J. of Algebra 223 (2000) 562--583.
}
:
does there exist a group which is
finitely presented, residually finite, and perfect?
\newline
Not clear we should keep this !!!!!!!!!!!!!!!!!!!!

\subhead
Un vieux bout de (7.B)
\endsubhead

R\'esultat d'une tentative non concluante avec J.W~;
\`a ne pas conserver ici tel quel
!!!!!!!!!!!!!!!!!!!!!!!!!!
\newline
For this footnote, write $\overline{A}$ instead of $P_{\hat 2}(A)$.
Let $N$ be a normal subgroup of $\Delta$, not reduced to one element.
Then $\overline{N}$ is a normal subgroup of $\overline{\Delta} =  
\overline{\Gamma}$,
not reduced to one element,
so that the quotient $\overline{\Delta}/\overline{N}$
is a finite $2$-group
(Theorem 12.3.31 in \cite{LeGMc--02}).
In particular, the number of cosets of $\overline{\Delta}$ modulo $ 
\overline{N}$
is finite; hence $\overline{N}\Delta$ is closed in $\overline{\Delta}$,
and $\overline{N}\Delta = \overline{\Delta}$.
It follows from the \lq\lq second isomorphism theorem\rq\rq \ that
$\overline{\Delta}/\overline{N}$ is isomorphic to
$\Delta/(\overline{N} \cap \Delta)$.
\newline
??????? S'il \'etait vrai que
$\overline{N} \cap \Delta = N$,
nous aurions montr\'e que tout quotient propre de $\Delta$
est un quotient propre de $\overline{\Delta} = \overline{\Gamma}$,
en particulier que tout quotient propre de $\Delta$ est un $2$-groupe  
fini.

\subhead
(7.C)
\endsubhead

R\'ef. \cite{Stall--65}~: v\'erif. intitul\'e. Par ailleurs, th.~3.4,  
apparemment~:

Soit $\varphi : A \longrightarrow B$ un homomorphisme de groupes
qui induit un iso sur $H_1(-,\bold Z)$ et un epi sur $H_2(-,\bold Z)$~;
alors $\varphi$ induit un iso $A/C^n(A) \longrightarrow B/C^n(B)$  
pour tout $n$.

\subhead
A informer une fois de notre travail
\endsubhead

cochran\@math.rice.edu et sharvey\@math.mit.edu
Voir

T. Cochran and S. Harvey, {\it Homology of derived series of groups,}

arXiv:math.GT/0407203 v1  12 Jul 2004.

\bigskip
\head{\bf
Reminder for GA + DK + PH -- End of July, 2005 (???)
}\endhead

Check (2.D) once more.

\medskip

Fix diagrams in (3.A) and (5.C).

\medskip

For (6.B), see wether we can find in the literature a finitely presented
example (namely a non-Hopfian f.p. group which is soluble of class $3$).

\medskip

Example of free products of finite soluble groups for (6.C) ?
More precisely, given four finite soluble groups $A,B,C,D$,
decide when $P \Cal S (A \ast B)$ and $P \Cal S (C \ast D)$
are isomorphic.
(Recall that these groups are residually soluble,
by the result of Gruenberg recalled as a footnote to (6.E).)

\medskip

In (6.D), there used to be the following remark/question, as a footnote.

Modulo checking: for $\Gamma_d(p)$, the true prosoluble completion,
the prosoluble completion, and the profinite completion all coincide
with the {\it pro-$p$-completion} $P  \hat p  (\Gamma)$.

\medskip

Understand other examples of true prosoluble kernels,
e.g. those listed in  (7.F).

\medskip

It seems true that P. Hall has introduced the WORDS \lq\lq residually  
P\rq\rq ,
but the notion has apparently been used earlier by others,
including R. Baer (Bull. AMS 50, 1944, 143--160)
and Mal'cev \cite{Malc'--49}.

\medskip

??? Anything in Goulnara's file from early May which is not in the  
present one ???

\medskip

About (7.C) and the definition of \lq\lq parafree groups\rq\rq .
Say a group $\Gamma$ has properties
\roster
\item"(P1)" if there  exists a free group $F$ such that $\Gamma/C^j 
\Gamma$ and
$F/C^jF$ are  isomorphic for all~$j$,
\item"(P2)" if there exists a group homomorphism from a free group $F 
$ to $\Gamma$
which  induces isomorphisms between  $F/C^jF$ and $\Gamma/C^j\Gamma$  
for all $j$.
\endroster
The question is to know if Property (P1) implies Property (P2).

\smallskip

The question has been asked to Gilbert Baumslag (mail of August 2).
In case he answers, change (7.C) accordingly.

\medskip

About question in (7.C), study the basilica group. Reminder and one  
precise question
below.

\medskip
\vskip.2cm
\Refs
\nofrills{}
\widestnumber\no{MaKaS--66}

Denote by $T$ the rooted binary tree.
Its vertex set is the set $\{0,1\}^*$ of finite words spelled with $0 
$ and $1$'s;
each vertex $j = (j_1,\hdots,j_n) \in \{0,1\}^*$ is connected by an edge
to the two vertices $0j = (0,j_1,\hdots,j_n)$ and $1j$,
as well as to $(j_2,\hdots,j_n)$ if $j$ is not the root vertex  
(namely the empty word).
The {\it flip} $\sigma$ is the automorphism of $T$ defined by
$\sigma(0j) = 1j$ and $\sigma(1j) = 0j$. Define recursively
\footnote{
Check if the notation here coincides with those of L.B., who writes
$$
a = <<b,1>>(1,2) \qquad \text{and} \qquad b = <<a,1>>,
$$
and V.N., who uses similar notation,
while \cite{GrZu--02a} use $a = (1,b)$ and $b = (1,a)\epsilon$.
}
two automorphisms $a,b$ of $T$ by
$$
\matrix
&\operatorname{a}(0j) \, &= \, &1j
        &\hskip1cm &\operatorname{b}(0j) \, &= \, &0\operatorname{a} 
(j) \\
&\operatorname{a}(1j) \, &= \, &0b(j)
         &\hskip1cm &\operatorname{b}(1j) \, &= \, &1j
\endmatrix
$$
or, in shorter and by now classical notation, by
$$
a = (b,1)\sigma \qquad \text{and} \qquad b \, = \, (a,1).
$$
The {\it basilica group} is the group of automorphisms of $T$  
generated by $a$ and $b$;
it is denoted by $\operatorname{IMG}(z^2-1)$, where IMG stands for
\lq\lq iterated monodromy group\rq\rq , a terminology explained in  
\cite{Nekra}.
The following properties have been established.
\roster
\item"(i)" The group $\operatorname{IMG}(z^2-1)$ is torsion free.
\item"(ii)" The semi-group generated by $a$ and $b$ is free; in  
particular
    $\operatorname{IMG}(z^2-1)$ is of exponential growth.
\item"(iii)" The group $\operatorname{IMG}(z^2-1)$ does not contain
   non-abelian free subgroups.
\item"(iv)" The Schur multiplier $H_2(\operatorname{IMG}(z^2-1),\bold  
Z)$
    is free abelian of infinite rank; in parti\-cular, $\operatorname 
{IMG}(z^2-1)$ is not
    finitely presented.
\item"(v)" The group $\operatorname{IMG}(z^2-1)$ is just non soluble,  
namely is non
     soluble with all proper quotients soluble;
    more precisely, any proper quotient of this group is an extension  
of a nilpotent group
    by a finite $2$-group,
\item"(vi)" The group $\operatorname{IMG}(z^2-1)$ is amenable.
\endroster
See \cite{GrZu--02a}, \cite{GrZu--02b}, \cite{BarVi}, \cite 
{BaGrB--03}, and \cite{BarGr}.
(And, possibly, also the book \cite{GrZuk}.)

Set $\Gamma = \operatorname{IMG}(z^2-1)$.
Observe that  $\Gamma$ is residually finite,
since the automorphism group of $T$ is residually finite;
thus $\Gamma$ embeds in $P \Cal F (\Gamma)$.
It follows from Property (v) above that
$\Gamma$ embeds in both $P \Cal S (\Gamma)$ and $P \Cal F \Cal S  
(\Gamma)$,
and also that $P \Cal F (\Gamma) = P \Cal F \Cal S (\Gamma)$.
Recall from (4.B) that there is a canonical homomorphism
$P \Cal S (\Gamma) \longrightarrow P \Cal F \Cal S (\Gamma)$.

{\bf Question:}
is $P \Cal S (\Gamma) \longrightarrow P \Cal F \Cal S (\Gamma)$
an isomorphism? if not, what are its kernel and cokernel?

First question for us: identify an infinite soluble quotient of $ 
\Gamma$.

\bigskip

\ref \no BaGrN--03 \by L. Bartholdi, R. Grigorchuk, and V. Nekrashevych
\paper From fractal groups to fractal sets
\jour in \lq\lq Fractals in Graz\rq\rq
\yr Birkh\"auser, 2003 \pages 25--118
\endref

\ref \no BarGr \by L. Bartholdi and R. Grigorchuk
\paper On a group associated to $z^2-1$
\jour Preprint
\endref

\ref \no BarVi \by L. Bartholdi and B. Virag
\paper Amenability via random walks
\endref

\ref \no GrZu--02a \by R. Grigorchuk and A. Zuk
\paper On a torsion-free weakly branch group defined by a three state  
automaton
\jour Internat. J. Algebra Comput \vol 12 \yr 2002 \pages 223--246
\endref

\ref \no GrZu--02b \by R. Grigorchuk and A. Zuk
\paper On the sprectrum of a torsion-free weakly branch group defined  
by a three state
automaton
\jour in \lq\lq Computational and Statistical Group Theory, Contemp.  
Math.
\vol ??? \yr Amer. Math. Soc. ???
\endref

\ref \no GrZuk \by R. Grigorchuk and A. Zuk
\book Automata groups, their spectra and classification
\publ Book in preparation
\endref

\ref \no Nekra \by V. Nekrashevych
\book Iterated monodromy groups
\yr book to appear in 2005
\endref

\endRefs

\vskip1cm

\subhead
Quelques questions
\endsubhead

   Voici encore des exercices pour le prochain jour de pluie.
A ne pas mettre pour l'instant dans la redaction !!!

\medskip

Soit P l'une des proprietes suivantes.
Est-ce que la propriete P passe d'un groupe G
(suppose residuellement resoluble pour eviter des absurdites)
a sa completion resoluble pleine PS(G) ?

\medskip

Exemples de P :

tout sous-groupe resoluble est abelien,

tout sous-groupe resoluble est virtuellement abelien,

tout sous-groupe resoluble est cyclique,

tout sous-groupe resoluble est nilpotent,

tout sous-groupe resoluble est nilpotent de classe de nilpotence au  
plus k,

tout sous-groupe resoluble est resoluble de classe de solubilite au  
plus k,

etc, etc, etc, etc, etc, etc, etc, etc, etc, etc, etc, etc, etc, etc,  
etc, etc,
etc, etc, etc, etc, etc, etc.

\bigskip

\subhead
Torsion
\endsubhead

Is it true that $P \Cal S (G)$ is never a torsion group (excluding  
the case of a
finite group $G$~!). For $P \Cal F (G)$, see Zelmanov, {\it On  
compact periodic
groups}.

Remark. If $G$ is torsion free, then $P \Cal F (G)$ has torsion (...).
See Lubotzky, {\it Torsion in profinite completions of torsion-free  
groups.}
See also J.R. McMullen, {\it The profinite completion of certain  
residually
finite-$p$-groups}, Monatsheft Math. 99 (1985) 311--314.

\subhead
Sous-groupe parfait maximal
\endsubhead

    (5.F) Let $\Gamma$ be a countable  group.
If there exists in $\Gamma$ a perfect subgroup not reduced to $\{1\}$,
then clearly $\Gamma$ is not soluble
(indeed not residually soluble).
Recall however that the converse {\it does not hold!}
in particular, $\Gamma$ can be non-soluble,
namely $D^{\infty}(\Gamma) = \cap_{n  \ge 1}D^n(\Gamma)$
can be not reduced to~$\{1\}$,
and yet $D^{\infty}(\Gamma)$ can be soluble.
More generally, Mal'cev has constructed for each
???countable??? ordinal $\alpha$ a countable group $\Gamma$
such that $D^{\alpha}(\Gamma) = \{1\}$,
and $D^{\beta}(\Gamma) \ne \{1\}$ whenever $\beta < \alpha$
\cite{?????Mal'cev?????}.


\bigskip
\Refs
\widestnumber\no{BruBW--79}

\ref \no Abels--79 \by H. Abels
\paper An example of a finitely presented solvable group
\jour in \lq\lq Homological group theory\rq\rq , Durham 1977,
C.T.C. Wall Editor \yr Cambridge Univ. Press 1979 \pages 205--211
\endref

\ref \no Baums--67 \by G.~Baumslag
\paper Some groups that are just free
\jour Bull. Amer. Math. Soc. \vol 73 \yr 1967 \pages 621--622
\endref

\ref \no Baums--68 \by G.~Baumslag
\paper More groups that just about free
\jour Bull. Amer. Math. Soc. \vol 74 \yr 1968 \pages 752--754
\endref

\ref \no Baums--69 \by G.~Baumslag
\paper A non--cyclic one--relator group
all of whose finite quotients are cyclic
\jour J. Austral. Math. Soc. \vol 10 \yr 1969 \pages 497--498
\endref

\ref \no Baums--71 \by G.~Baumslag
\paper Positive one--relator groups
\jour Trans. Amer. Math. Soc. \vol 156 \yr 1971 \pages 165--183
\endref


\ref \no Baums--05 \by G.~Baumslag
\paper Parafree groups
\inbook Infinite Groups: Geometric, Combinatorial and Dynamical Aspects
\eds L.~Bartholdi, T.~Ceccherini-Silberstein, T.~Smirnova-Nagnibeda, and A.~Zuk
\moreref Series: Progress in Mathematics 
\vol 248
\yr 2005  \pages 1-15
\publ Birkh\"auser
\endref


\ref \no Birkh--37 \by G. Birkhoff
\paper More--Smith convergence in general topology
\jour Ann. of Math. \vol 38 \yr 1937 \pages 39--56
\endref

\ref \no Bourb--60 \by N.~Bourbaki
\book Topologie g\'en\'erale, chapitres 3 et 4, troisi\`eme \'edition
\publ Hermann \yr 1960
\endref

\ref \no Bourb--70 \by N.~Bourbaki
\book Alg\`ebre, chapitres 1 \`a 3
\publ Diffusion C.C.L.S., Paris \yr 1970
\endref

\ref \no Bourb--72 \by N.~Bourbaki
\book Groupes et alg\`ebres de Lie, chapitres 2 et 3
\publ Hermann \yr 1972
\endref

\ref \no Bourb--82 \by N.~Bourbaki
\book Groupes et alg\`ebres de Lie, chapitre 9
\publ Masson \yr 1982
\endref

\ref \no BriGr--04 \by M.~Bridson and F.J.~Grunewald
\paper Grothendieck's problems concerning profinite completions and
representations of groups
\jour Annals of Math.  \vol 160 \yr 2004 \pages 359--373
\endref

\ref \no BruBW--79 \by A.M. Brunner, R.G. Burns, and J. Wiegold
\paper On the number of quotients, of one way or another,
of the modular group
\jour Math. Scientist \vol 4 \yr 1979 \pages 93--98
\endref


\ref \no BruSV--99 \by A.M. Brunner, S. Sidki, and A.C. Vieira
\paper A just nonsolvable torsion-free group defined on the binary tree
\jour J. of Algebra \vol 211 \yr 1999 \pages 99--114
\endref


\ref \no Clair--99 \by B. Clair
\paper Residual amenability and the approximation of $L^2$--invariants
\jour Michigan Math. J. \vol 46 \yr 1999 \pages 331--346
\endref

\ref \no CocHa \by T. Cochran and S. Harvey
\paper Homology and derived series of groups
\jour arXiv:math.GT/0407203
\endref

\ref \no DiSMS--91 \by J.D.~Dixon, M.~du Sautoy, A.~Mann, and D.~Segal
\book Analytic pro--$p$ groups
\publ Cambridge University Press \yr 1991
\endref

\ref \no EleSz \by G.~Elek and E.~Szabo
\paper On sofic groups
\jour arXiv:math.GR/0305352
\endref

\ref \no GlSoS \by Y.~Glasner, J.~Souto and P.~Storm
\paper Finitely generated subgroups of lattices in $PSL_2\bold C$
\jour arXiv:math.GT/0504441
\endref

\ref \no Grigo--80 	\by R.I. Grigorchuk
\paper Burnside's problem on periodic groups
\jour Functional Anal. Appl. \vol 14 \yr 1980 \pages 41--43
\endref



\ref \no Grigo--00	\by R.I. Grigorchuk
\paper Just infinite branch groups
\jour in [SaSeS--00]  \yr 2000 \pages 121--179
\endref

\ref \no Groth--70 \by A.~Grothendieck
\paper Repr\'esentations lin\'eaires
et compactification profinie des groupes discrets
\jour Manuscripta Math. \vol 2
\yr 1970 \pages 375--396
\endref

\ref \no Gruen--57 \by L. Gruenberg
\paper Residual properties of infinite soluble groups
\jour Proc. London Math. Soc. \vol 7 \yr 1957 \pages 29--62
\endref

\ref \no HallM--50 \by M.~Hall
\paper A topology for free groups and related groups
\jour Annals of Math. \vol 52 \yr 1950 \pages 127--139
\endref



\ref \no HallP--59 \by P.~Hall
\paper On the finiteness of certain soluble groups
\jour Proc. London Math. Soc. \vol 9 \yr 1959 \pages 595--622
[= Collected Works, 515--544]
\endref

\ref \no HallP--61 \by P.~Hall
\paper The Frattini subgroups of finitely generated groups
\jour Proc. London Math. Soc. \vol 11 \yr 1961 \pages 327--352
[= Collected Works, 581--608]
\endref

\ref \no Harpe--00 \by P. de la Harpe
\book Topics in geometric group theory
\publ The University of Chicago Press \yr 2000
\endref

\ref \no HigSt--54 \by G.~Higman and A.H.~Stone
\paper On inverse systems with trivial limits
\jour J. London Math. Soc. \vol 29 \pages 233--236
\endref

\ref \no Kahro1 \by D.~Kahrobaei
\paper On the residual solvability of generalized free products of
finitely
generated nilpotent groups
\jour arXiv:math.GR/0510465
\endref

\ref \no Kahro2 \by D.~Kahrobaei
\paper Are doubles of residually solvable groups, residually solvable?
\jour Preprint
\endref


\ref \no KasNi--06 \by  M.~Kassabov and N.~Nikolov
\paper  Cartesian products as profinite completions
\jour arXiv:math.GR/ 0602446
\endref


\ref \no Kelle--55 \by J.L.~Kelley
\book General topology
\publ Van Nostrand \yr 1955
\endref



\ref \no MacDu--46 \by C.C.~MacDuffee
\book The theory of matrices
\publ Chelsea Publ. Comp. \yr 1946
\endref

\ref \no MaKaS--66 \by W.~Magnus, A.~Karras, and D.~Solitar
\book Combinatorial group theory
\publ J. Wiley \yr 1966
\endref

\ref \no Mal'c--40 \by A.I. Mal'cev
\paper On the faithful representation of infinite groups by matrices
\jour Amer. Math. Soc. Transl (2) \vol 45 \yr 1965 \pages 1--18
[Russian original: Mat. SS.(N.S.) 8(50) (1940), pp.  405--422]
\endref

\ref \no Malc'--49 \by A.I. Mal'cev
\paper Generalized nilpotent algebras and their associated groups
\jour Mat. Sbornik N.S. \vol 25(67) \yr 1949 \pages 347--366
\endref

\ref \no NikSe--03 \by N. Nikolov and D. Segal
\paper Finite index subgroups in profinite groups
\jour C.R. Acad. Sci. Paris, S\'er. I \vol 337 \yr 2003 \pages 303--308
\endref

\ref \no NikSe--06a \by N. Nikolov and D. Segal
\paper On finitely generated profinite groups I: strong completeness and uniform bounds
\jour arXiv:math.GR/0604399
\endref 

\ref \no NikSe--06b \by N. Nikolov and D. Segal
\paper On finitely generated profinite groups II, products in quasisimple groups
\jour arXiv:math.GR/0604400
\endref



\ref \no Pet--73 \by H.L. Peterson
\paper Discontinuous characters and subgroups of finite index.
\jour Pacific J. Math. \vol 44 \yr 1973 \pages 683--691
\endref


\ref \no PlaTa--86 \by V.P.~Platonov and O.I.~Tavgen
\paper On Grothendieck's problem of profinite completions of groups
\jour Soviet Math. Dokl. \vol 33 \yr 1986 \pages 822--825
\endref

\ref \no PlaTa--90 \by V.P.~Platonov and O.I.~Tavgen
\paper Grothendieck's problem on profinite completions and  
representations of groups
\jour K--theory \vol 4 \yr 1990 \pages 89--101
\endref

\ref \no RapVa--89 \by E. Raptis and D. Varsos
\paper Residual properties of HNN--extensions with base group an abelian
group
\jour J. Pure Appl. Algebra \vol 59 \yr 1989 \pages 285--290
\endref

\ref \no RibZa--00 \by L.~Ribes and P.~Zalesskii
\book Profinite groups
\publ Springer \yr 2000
\endref

\ref \no Robin--82 \by D.J.S.~Robinson
\book A course in the theory of groups
\publ Springer \yr 1982
\endref

\ref \no SaSeS--00 \eds M. du Sautoy, D. Segal, and A. Shalev
\book New horizons in pro-$p$-groups
\publ Birkh\"auser \yr 2000
\endref


\ref \no Se--73/94 \by J--P.~Serre
\book Cohomologie galoisienne
\publ Lecture Notes in Math. 5, Springer \yr 1973
[cinqui\`eme \'edition 1994]
\endref

\ref \no Stall--65 \by J. Stallings
\paper Homology of central series of groups
\jour J. of Algebra \vol 2 \yr 1965 \pages 170--181
\endref

\ref \no Weil--40 \by A.~Weil
\book L'int\'egration dans les groupes topologiques et ses applications
\publ Hermann \yr 1940
\endref

\ref \no Wilso--98 \by J.~Wilson
\book Profinite groups
\publ Clarendon Press \yr 1998
\endref


\endRefs

\Refs\nofrills{}
\widestnumber\no{MaKaS--66}

\ref \no HallP--54 \by P.~Hall
\paper The splitting properties of relatively free groups
\jour Proc. London Math. Soc. \vol 4 \yr 1954 \pages 343--356
\endref
\endRefs

\vskip1cm

\Refs
\widestnumber\no{MaKaS--66}

\ref \no CocHa \by T. Cochran and S. Harvey
\paper Homology and derived series of groups
\jour arXiv:math GT/0407203/ 2004
\endref

\ref \no Hall \by P.~Hall
\paper The Frattini subgroups of finitely generated groups
\jour Proc. London Math. Soc. \vol 11 \yr 1961 \pages 327--352
[= Collected Works, 581--608]
\endref

\ref \no LewLi \by R.H. Lewis and S. Liriano
\paper Isomorphism classes and derived series of certain almost free  
groups
\jour Experimental Math. \vol 3:3 \yr 1994  \pages 255--288
\endref


\endRefs
\enddocument